\documentclass[12pt]{article}

\usepackage[utf8]{inputenc}
\usepackage{verbatim}
\usepackage{graphicx}
\usepackage{listings}
\usepackage{xpatch}

\usepackage{url}
\usepackage{hyperref}

\usepackage{pslatex}

\usepackage{datetime}
\newdateformat{versiondate}{%
\THEMONTH\THEDAY}

\usepackage{amsthm}
\usepackage{amsmath}
\usepackage{amssymb}
\usepackage{amsfonts}
\usepackage{authblk}

\newtheoremstyle{zoltanstyle}
  {1em} 
  {\topsep} 
  {} 
  {} 
  {\bfseries} 
  {.} 
  {.5em} 
  {} 

\theoremstyle{zoltanstyle}
\swapnumbers
\xpatchcmd\swappedhead{~}{.~}{}{}
\newtheorem{body}{}
\numberwithin{body}{section}
\newtheorem{definition}[body]{Definition}
\newtheorem{theorem}[body]{Theorem}
\newtheorem{proposition}[body]{Proposition}
\newtheorem{corollary}[body]{Corollary}

\newtheorem{lemma}[body]{Lemma}

\newtheorem{question}[body]{Question}

\expandafter\let\expandafter\oldproof\csname\string\proof\endcsname
\let\oldendproof\endproof
\renewenvironment{proof}[1][\proofname]{%
  \oldproof[\normalfont \bfseries #1.]%
}{\oldendproof}

\newcommand{\SetComp}[2]{\left\{ {#1}\:\middle|\:{#2} \right\}}   

\begin{document}

\title{Degree of satisfiability of\\ some special equations}
\author{Zoltan A. Kocsis\thanks{Department of Mathematics, The University of Manchester, Manchester UK}}
\date{\small 17 February 2020}

\maketitle

\begin{abstract}
A well-known theorem of Gustafson states that in a non-Abelian group the degree of satisfiability of $xy=yx$, i.e.~the probability that two uniformly randomly chosen group elements $x,y$ obey the equation $xy=yx$, is no larger than $\frac{5}{8}$. 

The seminal work of Antol\'{i}n, Martino and Ventura on generalizing the degree of satisfiability to finitely generated groups led to renewed interest in Gustafson-style properties of other equations. Positive results have recently been obtained for the 2-Engel and metabelian identities.

Here we show that the degree of satisfiability of the equations $xy^2=y^2x$, $xy^3=y^3x$ and $xy=yx^{-1}$ is either 1, or no larger than $1-\varepsilon$ for some positive constant $\varepsilon$. Using the Antol\'{i}n-Martino-Ventura formalism, we introduce criteria to identify which equations hold in a finite index subgroup precisely if they have positive degree of satisfiability. We deduce that the equations $xy=yx^{-1}$ and $xy^2=y^2x$ do not have this property.
\end{abstract}

\section{Introduction}

\begin{body}
In a finite group $G$, the probability that two uniformly randomly chosen elements of $G$ commute gives us a measure of how far $G$ is from being Abelian.
\end{body}

\begin{theorem}[Gustafson~\cite{gustafson-degree-of-commutativity}]\label{thm:gustafson-degree-of-commutativity}
In a non-Abelian finite group $G$, the probability that two uniformly randomly chosen elements of $G$ commute, i.e. $$\mathrm{dc}(G)  = \frac{ \left|\SetComp{(x,y) \in G^2}{xy=yx}\right| }{ \left|G^2\right| }$$
satisfies $\mathrm{dc}(G) \leq \frac{5}{8}$.
\end{theorem}

\begin{body}
Abelian groups have $\mathrm{dc}(G) = 1$, so Theorem~\ref{thm:gustafson-degree-of-commutativity} guarantees that there is no finite group with $\mathrm{dc}(G)$ strictly between $\frac{5}{8}$ and $1$. The bound $\frac{5}{8}$ is tight: the reader can check that the quaternion group $Q_8$ satisfies $\mathrm{dc}(Q_8) = \frac{5}{8}$.
\end{body}

\begin{definition}\label{def:degree-of-satisfiability}
Consider an equation $\varphi$ in the first-order language of group theory with $n$ free variables. We define the \textit{degree of satisfiability} of $\varphi$ as the quantity $$ \mathrm{ds}_G(\varphi) = \frac{ \left|\SetComp{\overline{x} \in G^n}{\varphi(\overline{x})}\right| }{ \left|G^n\right| }. $$
\end{definition}

\begin{body}
Clearly $\mathrm{dc}(G) = \mathrm{ds}_G(xy=yx)$. Many results are known about the distribution of $\mathrm{dc}(G)$ as $G$ ranges over the finite groups. In particular the groups $G$ satisfying $\mathrm{dc}(G) > \frac{11}{32}$ have been completely classified~\cite{rusin-classification}. However, the case of $\mathrm{ds}_G(\varphi)$ for an arbitrary equation $\varphi$ is much less well-understood. Even a general counterpart to Theorem~\ref{thm:gustafson-degree-of-commutativity} remains elusive.
\end{body}

\begin{question}\label{question:satisfiability-gap-for-ds}
Consider an equation $\varphi$ in the language of group theory. Can we find a constant $\varepsilon > 0$ such that for every group $G$, we have either $\mathrm{ds}_G(\varphi) = 1$ or else $\mathrm{ds}_G(\varphi) \leq 1 - \varepsilon < 1$?
\end{question}

\begin{body}
Delizia,~Jezernik,~Moravec~and~Nikotera~\cite{delizia-satisfiability} have obtained a positive answer to Question~\ref{question:satisfiability-gap-for-ds} in a special case when $\varphi$ is either $[x,[y,z]] = 1$ or $[[x,y],[z,w]] = 1$. Lescot~\cite{lescot-nilpotency} generalized Theorem~\ref{thm:gustafson-degree-of-commutativity} to the nested simple commutator equation $[x_1,x_2,\dots,x_n] = 1$, giving the constant $\varepsilon = \frac{3}{2^{n+2}}$.
\end{body}

\begin{body}
If Question~\ref{question:satisfiability-gap-for-ds} has a positive answer for some equation $\varphi$, we say that $\varphi$ \textit{has finite satisfiability gap}. The existence of an equation without a finite satisfiability gap is an open problem, even in the case where $\varphi$ has only one free variable. The case of $x^2 = 1$ can be settled by an elementary argument (Proposition~\ref{prop:square-bound}). The case of $x^3 = 1$ follows from a theorem of Laffey~\cite{laffey-powers-of-three}. Partial results are known for $x^p = 1$ with $p>3$. \cite{mann-survey}
\end{body}

\begin{proposition}\label{prop:square-bound}
The equation $x^2=1$ has finite satisfiability gap $\varepsilon = \frac{1}{4}$.
\begin{proof}
Assume that $\mathrm{ds}_G(x^2=1) \geq \frac{3}{4}$. We claim that all $x\in G$ with $x^2 = 1$ belong to the center. It suffices to prove that the centralizer of any such element is the entire group. Choose a random $y \in G$. We have $xy=yx$ if $y^2 = (xy)^2$. But the latter happens with probability $> \frac{1}{2}$, since our assumption guarantees that both the left hand and the right hand sides equal $1$ with probability $> \frac{3}{4}$. Thus, the centralizer of $x$ is larger than half of the group; by Lagrange's theorem, it must be the whole group. Therefore, all $x\in G$ with $x^2 = 1$ belong to the center. This means that the group $G$ is Abelian. The set $\SetComp{x \in G}{x^2 = 1}$ then forms a subgroup of $G$; by Lagrange's theorem, this subgroup is the whole group, so $\mathrm{ds}_G(x^2=1) = 1$.
\end{proof}
\end{proposition}

\begin{theorem}[Laffey~\cite{laffey-powers-of-three}]\label{thm:laffey-cube-bound}
Every finite group $G$ satisfies one of the following:
\begin{enumerate}
    \item Every element of $G$ has order 3,
    \item $G$ is a $3$-group and $\mathrm{ds}_G(x^3=1) \leq \frac{7}{9}$, or
    \item $G$ is not a $3$-group, and $\mathrm{ds}_G(x^3=1) \leq \frac{3}{4}$.
\end{enumerate}
Moreover, both of the bounds are tight.
\end{theorem}

\begin{body}
Gustafson~\cite{gustafson-degree-of-commutativity} proved that Theorem~\ref{thm:gustafson-degree-of-commutativity} holds in compact groups equipped with a left Haar measure. Much later, Antol\'{i}n, Martino and Ventura~\cite{antolin-infinite-dc} suggested a generalization of the degree of satisfiability to countably infinite groups, using sequences of densities (finitely additive measures) that ``measure index uniformly''\footnote{The precise definitions drifted over the years. We follow \cite{tointon-nilpotence} unless otherwise noted.}. They extended Theorem~\ref{thm:gustafson-degree-of-commutativity} to this new setting, and showed that in the virtually nilpotent case, the value of $\mathrm{dc}(G)$ is closely related to the algebraic structure of $G$ (Theorem~\ref{thm:antolin-commutativity}). 
\end{body}

\begin{definition}
Consider a group $G$ and a sequence $n \mapsto \mu_n$ of densities (finitely additive probability measures) on $G$. We say that $\mu_n$ \textit{measures index} if for any element $g \in G$ and finite-index subgroup $H < G$ we have $$\lim_{n \rightarrow \infty} \mu_n(gH) = \frac{1}{[G:H]}. $$

We say that $\mu_n$ \textit{measures index uniformly} if $\mu_n(gH)$ converges to $\frac{1}{[G:H]}$ uniformly for all elements $g \in G$ and subgroups $H < G$.
\end{definition}

\begin{body}
Densities that measure index uniformly encompass all the densities that have known infinite analogues of Theorem~\ref{thm:gustafson-degree-of-commutativity}: in particular, uniform measures of balls, simple random walks, and the Cesaro density~\cite{borovik-multiplicative-measures} all turn out to satisfy this definition. Moreover, it's an easy observation that pointwise Cartesian products of densities that measure index uniformly also measure index uniformly.
\end{body}

\begin{definition}
Take an equation $\varphi$ with $n$ free variables, and a group $G$ equipped with a sequence $i \mapsto \mu_i$ that measures index uniformly. We define the \textit{degree of satisfiability} of $\varphi$ with respect to the sequence $\mu$ as $$ \mathrm{ds}_{(G,\mu)}(\varphi) = \limsup_{n \rightarrow \infty} \mu_n^k \SetComp{ \overline{x} \in G^k }{ \varphi(\overline{x}) }.$$ 
\end{definition}

\begin{theorem}[Antol\'{i}n, Martino, Ventura~\cite{antolin-infinite-dc}]\label{thm:antolin-commutativity}
In a virtually nilpotent group $G$ equipped with a sequence $\mu$ that measures index uniformly, we have one of the following:
\begin{enumerate}
    \item $\mathrm{dc}(G,\mu) > 0$ and $xy=yx$ holds virtually in $G$, or
    \item $\mathrm{dc}(G,\mu) = 0$ and $xy=yx$ does not hold virtually in $G$.
\end{enumerate}
\end{theorem}

\begin{body}\label{body:antolin-discussion}
Antol\'{i}n, Martino and Ventura~\cite{antolin-infinite-dc} put forward that analogues of Theorem~\ref{thm:antolin-commutativity} should hold for most other equations. The suggested ``naive'' generalization (that $\mathrm{ds}_{(G,\mu)}(\varphi) > 0$ should hold precisely if $\varphi$ is a virtual law in $G$), clearly fails. Consider the equation $x^2 = 1$ in the infinite dihedral group $D_\infty$: the Abelian subgroups of $D_\infty$ are torsion-free, but their cosets satisfy the identity. Since the density cannot distinguish between cosets and subgroups, we will have $\mathrm{ds}_{(D_\infty,\mu)}(x^2=1) > 0$ for any $\mu$ that measures index uniformly. Fortunately, Martino, Tointon, Valiunas and Ventura~\cite{tointon-nilpotence} managed to show that this inability to distinguish between cosets and subgroups is the only obstacle in the way of generalizing Theorem~\ref{thm:antolin-commutativity} to arbitrary equations (Theorem~\ref{thm:tointon-big}).
\end{body}

\begin{theorem}[Martino, Tointon, Valiunas, Ventura~\cite{tointon-nilpotence}]\label{thm:tointon-big}
Consider an equation $\varphi$ in the language of group theory with $k$ free variables. Let $G$ be a finitely generated, virtually nilpotent group equipped with a sequence of densities $\mu$ that measures index uniformly. Then $\mathrm{ds}_{(G,\mu)}(\varphi) > 0$ precisely if there is a finite index subgroup $H<G$ and elements $g_1,\dots g_k$ satisfying $$\forall x_1 \in g_1 H. \dots \forall x_k \in g_k H. \varphi(x_1,\dots,x_k).$$
\end{theorem}

\begin{body}
In light of Theorem~\ref{thm:gustafson-degree-of-commutativity}, Proposition~\ref{prop:square-bound} and Theorem~\ref{thm:laffey-cube-bound}, it is natural to consider the equations $xy = y^{-1}x$, $xy^2 = y^2x$ and $xy^3 = y^3x$. In this paper we show that all three equations have finite satisfiability gap (Propositions~\ref{prop:xyy-equals-yyx-gap}, \ref{prop:xyyy-equals-yyyx-gap} and \ref{prop:xy-equals-yX-gap}). Our techniques yield a new proof of Lascot's bound on the satisfiability of the nested commutator equation as well (Corollary~\ref{cor:nested-commutator-gap}). We discuss conditions on equations $\varphi$ that satisfy not only Theorem~\ref{thm:tointon-big}, but a more direct analogue of Theorem~\ref{thm:antolin-commutativity} as well, in the sense that if $\varphi$ holds with positive probability, then it is in fact a virtual law.
\end{body}

\section{Results}\label{sec:results}

\begin{body}
We begin by determining a tight finite satisfiability gap for $xy^2=y^2x$. The workhorse of the proof is Lemma~\ref{lemma:normal-subgroup-bound}, an elementary observation about the relationship between the degree of satisfiability of $f(x) = 1$ and the probability that a normal subgroup contains elements of $\mathrm{im} f$.
\end{body}

\begin{lemma}\label{lemma:normal-subgroup-bound}
Consider a finite group $G$ and a function $f: G^n \rightarrow G$ for some $n \in \mathbb{N}$ such that the equation $f(x_1,x_2,\dots,x_n) = 1$ is first-order in the language of group theory, with satisfiability gap $\varepsilon$. For each normal subgroup $N \trianglelefteq G$, one of the following holds:
\begin{enumerate}
    \item $N$ contains the entire image of $f$, or
    \item for uniformly random group elements $x_1,\dots,x_n$, $f(x_1,\dots,x_n)$ belongs to $N$ with probability at most $1 - \varepsilon$.
\end{enumerate}
\begin{proof}
Assume that for uniformly randomly chosen group elements $x_1,\dots,x_n$, the value $f(x_1,\dots,x_n)$ belongs to $N$ with probability strictly larger than $1-\varepsilon$. Consider the quotient group $G/N$. Since each congruence class has the same size, taking the congruence class $xN$ of a uniformly random element $x \in G$ yields a uniformly random element of $G/N$. So pick $x_1, x_2, \dots, x_n$ uniformly randomly from $G$. Evaluate $f(x_1,x_2,\dots,x_n)$. With probability over $1-\varepsilon$, $f(x_1,x_2,\dots,x_n)$ belongs to $N$, and then $$f(x_1,x_2,\dots,x_n)N = f(x_1N,x_2N,\dots,x_nN) = N.$$
Consequently, the equation $f(x_1,x_2,\dots,x_n) = 1$ has degree of satisfiability larger than $1 - \varepsilon$ in $G/N$. Under our assumption on the satisfiability gap, we get that it has degree of satisfiability $1$. But then $\mathrm{im} f \subseteq N$.
\end{proof}
\end{lemma}

\begin{proposition}\label{prop:xphi-equals-phix-gap}
Consider a finite group $G$ and a function $f: G^n \rightarrow G$ for some $n \in \mathbb{N}$ such that the equation $f(x_1,x_2,\dots,x_n) = 1$ is first-order in the language of group theory, with satisfiability gap $\varepsilon$. If $x$ is a variable disjoint from $\{x_1,\dots,x_n\}$, then the equation $xf(x_1,\dots,x_n) = f(x_1,\dots,x_n)x$ has finite satisfiability gap~$\frac{1}{2}\varepsilon$.
\begin{proof}
Take a finite group $G$. If $f = f(x_1,\dots,x_n)$ belongs to the center $Z(G)$ for each $x_1,\dots,x_n \in G$, then $\mathrm{ds}_G(xf=fx)$ is 1, and we're done.  Otherwise, we can apply Lemma~\ref{lemma:normal-subgroup-bound} and the gap $\varepsilon > 0$ to conclude that for uniformly randomly chosen elements $x_1,\dots,x_n \in G$, the quantity $f(x_1,\dots,x_n)$ belongs to $Z(G)$ with probability at most $1 - \varepsilon$.

So pick $x,x_1,\dots,x_n \in G$ uniformly randomly. The event $f(x_1,\dots,x_n) \not\in Z(G)$ happens with probability at least $\varepsilon$. Moreover, given $f(x_1,\dots,x_n) \not\in Z(G)$, the centralizer $C = C_G\!\left(f(x_1,\dots,x_n)\right)$ exhausts at most half the group (\textit{dixit} Lagrange!), and the event $x \not\in C$ happens with probability at least~$\frac{1}{2}$. So both events happen with probability at least~$\frac{1}{2}\varepsilon$. But if both events happen, then $xf \neq fx$.
\end{proof}
\end{proposition}

\begin{corollary}\label{cor:nested-commutator-gap}
The nested simple commutator equation $\left[x_1,x_2,\dots,x_n\right] = 1$ has finite satisfiability gap $\varepsilon = \frac{3}{2^{n+2}}$.
\begin{proof}
Induction on $n$ using Theorem~\ref{thm:gustafson-degree-of-commutativity} and Proposition~\ref{prop:xphi-equals-phix-gap}.
\end{proof}
\end{corollary}

\begin{proposition}\label{prop:xyy-equals-yyx-gap}
The equation $xy^2 = y^2x$ has finite satisfiability gap $\varepsilon = \frac{1}{8}$.
\begin{proof}
Notice that by Proposition~\ref{prop:square-bound}, the equation $y^2 = 1$ has finite satisfiability gap $\frac{1}{4}$. Apply Proposition~\ref{prop:xphi-equals-phix-gap} to $f(y) = y^2$.
\end{proof}
\end{proposition}

\begin{body}
We now show the tightness of the gap bound determined in Proposition~\ref{prop:xyy-equals-yyx-gap} by constructing a group~$G$ with $\mathrm{ds}_G(xy^2=y^2x) = \frac{7}{8}$.
\end{body}

\begin{proposition}
The dihedral group $G = D_{16}$ has $\mathrm{ds}_G(xy^2=y^2x) = \frac{7}{8}$.
\begin{proof}
One could of course verify this exhaustively, but it's far quicker to consider the usual presentation $$\langle R, F \:|\: R^8 = 1, F^2 = 1, FR = R^{-1}F \rangle.$$
We immediately see that the squares are $1, R^2, R^4, R^6$ with $R^4$ generating the center. We know that $D_{16}/Z(D_{16}) \cong D_8$ which has $x^2 = 1$ with probability $\frac{3}{4}$ as desired. So we need to check the centralizers of the remaining squares, $R^2$ and $R^6$. It's immediate that all $R^n$ belong to these centralizers, but no element of the form $R^nF$ does, so they each constitute half the group.
\end{proof}
\end{proposition}

\begin{body}
Repeating the exact same argument as above, we would obtain the elegant, but not tight, gap $\varepsilon = \frac{1}{9}$ for the equation $xy^3 = y^3x$. Deriving a tight bound requires a further case analysis on the three possibilities of Theorem~\ref{thm:laffey-cube-bound}.
\end{body}

\begin{proposition}\label{prop:xyyy-equals-yyyx-gap}
The equation $xy^3 = y^3x$ has finite satisfiability gap $\varepsilon = \frac{4}{27}$.
\begin{proof}
Take a finite group $G$. As usual, if $x^3$ belongs to the center $Z(G)$ for each $x \in G$, then $ds_G(xy^3=y^3x)$ is 1. In the remaining cases, we proceed analogously to the proof of Proposition~\ref{prop:xyy-equals-yyx-gap}, by noticing that the event $xy^3 \neq y^3x$ happens precisely if both $y^3 \not\in Z(G)$ and $x \not\in C_G(y^3)$ happen. In accordance with Theorem~\ref{thm:laffey-cube-bound}, we need to choose $G$ as a $3$-group to minimize the probability of both events happening. Then by Lemma~\ref{lemma:normal-subgroup-bound}, the value $x^3$ belongs to $Z(G)$ with probability at most $\frac{7}{9}$.

Pick $x,y \in G$ uniformly randomly. The event $y^3 \not\in Z(G)$ happens with probability at least $\frac{2}{9}$. Moreover, since $G$ is a $3$-group and $y^3 \not\in Z(G)$, the centralizer $C_G(y^3)$ exhausts at most one third of the group, so the event $x \not\in C_G(y^3)$ happens with probability at least $\frac{2}{3}$. So these both happen with probability at least $\frac{4}{27}$.
\end{proof}
\end{proposition}

\begin{body}
Again, we wish to show the tightness of the gap bound derived in Proposition~\ref{prop:xyyy-equals-yyyx-gap}. Since the $\frac{2}{9}$ gap of Theorem~\ref{thm:laffey-cube-bound} is not achieved in any group of order less than 81, and the gap for $xy^3=y^3x$ can only be achieved for a $3$-group, our example will have to have order at least 243.
\end{body}

\begin{proposition}\label{prop:xyyy-equals-yyyx-tight}
There is a group $G$ of order 243 with $\mathrm{ds}_G(xy^3=y^3x) = \frac{23}{27}$.
\begin{proof}
Consider the finitely presented group $G$ on three generators $X,Y,F$, subject to the following relations:
\begin{enumerate}
    \item $X^9 = Y^9 = F^3 = 1$,
    \item $XY=YX$,
    \item $XF = FYX$, and
    \item $YF = FYX^6$.
\end{enumerate}
We can write each element of the group uniquely in the normal form $F^\ell Y^m X^n$ with natural number indices $\ell < 3$, $m < 9$ and $n < 9$. Denoting this element as $(\ell,m,n)$, we arrive at the following computation rule for the group operation:
\begin{align*}
        & (\ell_1,y_1,x_1) (0,y_2,x_2) = (\ell_1, y_2 + y_1, x_2 + x_1), \\
        & (\ell_1,y_1,x_1) (1,y_2,x_2) = (\ell_1 + 1, x_1 + y_1, x_1 + 6y_1)(0,y_2,x_2), \\
        & (\ell_1,y_1,x_1) (2,y_2,x_2) = (\ell_1 + 1, x_1 + y_1, x_1 + 6y_1)(1,y_2,x_2).
\end{align*}
where we perform the additions in the first component modulo 3, and in the other components modulo 9. This amounts to writing the group as a semidirect product of the form $(\mathbb{Z}/9\mathbb{Z})^2 \rtimes \mathbb{Z}/3\mathbb{Z}$. The explicit use of normal forms keeps the computations manageable.

Notice first the following: $(1,y,x)^3 = (0,3x+9y, 9x+18y) = (0,3x,0)$ and similarly for $(2,y,x)$. From the normal forms it follows immediately that $Y^3$ generates $Z(G)$, so the quotient $G/Z(G)$ indeed has 81 elements. Moreover, we get corresponding normal forms in the factor group (where one performs addition in the second component modulo 3): we denote these normal forms using square brackets, as $[\ell,y,x]$.

Thanks to the observations above, we know that elements of the form $[1,y,x]$ and $[2,y,x]$ always satisfy $[\ell,y,x]^3 = [0,0,0]$ (a la $FR^n$ in a dihedral group!). We have that $[0,y,x]^3 = [0,3y,3x]=[0,0,x]$, so an element of the form $[0,y,x]$ satisfies $x^3 = 1$ in $G/Z(G)$ precisely if $3x$ is a multiple of 9, i.e. $x \in \{0,3,6\}$. We conclude that the equation $x^3 = 1$ has degree of satisfiability $\frac{2}{3} + \frac{1}{9} = \frac{7}{9}$ in $G/Z(G)$.

The cubes of $G$ that do not belong to $Z(G)$ are the following:  $$(0,0,3),(0,0,6),(0,3,3),(0,3,6),(0,6,3),(0,6,6).$$
But these commute exactly with the elements of the form $(0,y,x)$, and therefore with probability $\frac{1}{3}$. Consequently, the probability that $xy^3 = y^3x$ is satisfied for uniformly randomly selected elements $x,y \in G$ is $1 - \left(1 - \frac{7}{9}\right)\left(1 - \frac{1}{3}\right) = \frac{23}{27}$. 
\end{proof}
\end{proposition}

\begin{body}
Finally, we turn to the equation $xy=yx^{-1}$. The proof of the finite satisfiability gap result is more elaborate in this case, combining features of previous ones: in the Abelian case, the bounds obtained in Proposition~\ref{prop:square-bound}, and in the non-Abelian case the bound of Theorem~\ref{thm:gustafson-degree-of-commutativity}. Unlike the previous cases, however, the tightness proof is delightfully simple: the bound is obtained in the dihedral group~$D_8$.
\end{body}

\begin{proposition}\label{prop:xy-equals-yX-gap}
The equation $xy=yx^{-1}$ has finite satisfiability gap $\varepsilon = \frac{3}{8}$.
\begin{proof}
First assume that $G$ is Abelian. Then $xy=yx^{-1} \leftrightarrow xy=x^{-1}y \leftrightarrow x^2 = 1$. But the set $\{x \in G \:|\: x^2 = 1\}$ forms a subgroup of the Abelian group $G$. By Lagrange's theorem we have either $\mathrm{ds}_G(x^2=1) = 1$ or $\mathrm{ds}_G(x^2=1) \leq \frac{1}{2}$, and in the latter case $\mathrm{ds}_G(xy=yx^{-1}) \leq \frac{1}{2} \leq \frac{5}{8}$ as desired.

For the remainder of this proof, assume that $G$ is not Abelian. Consider the set $ZR(G) =\{x \in G \:|\: \forall y \in G. xy = yx^{-1}\}$. Take arbitrary $a \in ZR(G)$. By definition we have $\forall y \in G. ay = ya^{-1}$. So we can take $y=a$ and obtain $a^2 = 1$. This means that in fact $\forall y \in G. ay=ya$, and thus $a$ belongs to the center $Z(G)$. We conclude that $ZR(G) \subseteq Z(G)$. Note that since $G$ is not Abelian, a random $x \in G$ belongs to $Z(G)$, and hence to $ZR(G)$ with probability at most $\frac{1}{4}$.

Now consider the set $CL_G(a)=\{y \in G \:|\: ay=ya^{-1} \}$. Generally, $CL_G(a)$ does not arise as a subset of the centralizer $C_G(a)$. However, if $CL_G(a) \neq \emptyset$ then we have $|CL_G(a)| = |C_G(a)|$. First of all, we observe that if $z \in C_G(a)$ and $y \in CL_G(a)$ then $yz \in CL_G(a)$. Conversely, if $y_1, y_2 \in CL_G(a)$, then $ay_1y_2^{-1} = y_1a^{-1}y_2^{-1} = y_1y_2^{-1}a$, so any two elements of $CL_G(a)$ differ by some element of the centralizer $C_G(a)$. With this, we establish $|CL_G(a)| = |C_G(a)|$. Since the group $G$ is not Abelian, if $ZR(a) \neq G$, then a random $x \in G$ belongs to $CL_G(a)$ with probability at most $\frac{1}{2}$.

Pick random $x,y \in G$. With probability less than $\frac{1}{4}$, we have $x \in ZR(G)$, and then $xy=yx^{-1}$, regardless of $y$. Otherwise (with probability at least $\frac{3}{4}$), $x \not\in ZR(G)$ and then $xy=yx^{-1}$ holds precisely if $y \in CL_G(x)$, which fails with probability at least $\frac{1}{2}$. Thus, $xy=yx^{-1}$ fails with probability at least $\frac{3}{4} \cdot \frac{1}{2} = \frac{3}{8}$.
\end{proof}
\end{proposition}

\begin{body}
As we noted in \ref{body:antolin-discussion}, density cannot distinguish a subgroups from its cosets: this is why identities that hold with positive probability are coset identities, but usually not virtual laws. One seeks a criterion that can distinguish between equations that satisfy direct analogues of Theorem~\ref{thm:antolin-commutativity} (i.e. equations which are virtual laws precisely if they hold with positive probability) from those that satisfy only the weaker conclusion of Theorem~\ref{thm:tointon-big}. Proposition~\ref{prop:commutative-transfer} provides a legible reason why the equation $xy=yx$ belongs to the former class. The equations $xy=yx^{-1}$ and $xy^2=y^2x$ do not pass the same muster, and indeed we can find groups in which the analogues of Theorem~\ref{thm:antolin-commutativity} fail for these equations.
\end{body}

\begin{proposition}\label{prop:commutative-transfer}
Consider a group $G$ with a finite index subgroup $H<G$ and elements $a,b \in G$ satisfying $\forall x \in aH. \forall y \in bH. xy=yx$. Then $G$ is virtually Abelian.
\begin{proof}
We prove that $xy=yx$ holds universally in $H$ itself. By assumption, we have $\forall x, y \in H. axby = byax$. Since $1 \in H$, we can set $x=y=1$ to obtain $ab=ba$. Now, obtain $aby=bya$ by setting $x=1$ in $axby=byax$. Using $ab=ba$, we get $bay=bya$ and hence $ay=ya$ for all $y \in H$. Similarly, obtain $axb=bax$ by setting $y=1$ in the original identity. Use $ab=ba$ to get $axb=abx$, and hence $xb=bx$ for all $x \in H$. Together these give $abxy = axby = aybx = abyx$, and thus $xy=yx$.
\end{proof}
\end{proposition}

\begin{body}
We obtain Theorem~\ref{thm:antolin-commutativity} as an immediate corollary of Theorem~\ref{thm:tointon-big} using our Proposition~\ref{prop:commutative-transfer}. The result motivates the following definition.
\end{body}

\begin{definition}\label{def:transfer}
Consider an equation $\varphi$ in the language of group theory with $k$ free variables. We say that $\varphi$ \textit{transfers} (from cosets to subgroups) if for each group $G$ and subgroup $H<G$ the following holds: if there are elements $g_1,\dots g_k$ satisfying $$\forall x_1 \in g_1 H. \dots \forall x_k \in g_k H. \varphi(x_1,\dots,x_k)$$
then $\varphi$ is a virtual law in $H$.
\end{definition}

\begin{body}
In light of Definition~\ref{def:transfer}, we can recast Proposition~\ref{prop:commutative-transfer} as ``the equation $xy=yx$ transfers''. It's clear that if an equation $\varphi$ transfers, then $\varphi$ satisfies the direct analogue of Theorem~\ref{thm:antolin-commutativity} stated below.
\end{body}

\begin{proposition}\label{prop:transfer-big}
Consider an equation $\varphi$ in the language of group theory with $k$ free variables. Assume that $\varphi$ transfers. Let $G$ be a finitely generated, virtually nilpotent group equipped with a sequence of densities $\mu$ that measures index uniformly. Then $\mathrm{ds}_{(G,\mu)}(\varphi) > 0$ precisely if $\varphi$ holds virtually in $G$.
\end{proposition}

\begin{body}
By considering small cyclic groups, the reader can quickly prove that the equation $xy=yx^{-1}$ does not transfer. The equation $xy^2 = y^2x$ is a tougher nut to crack: we now prove that it does not transfer by constructing an infinite virtually nilpotent group that does not satisfy the conclusion of Proposition~\ref{prop:transfer-big}.
\end{body}

\begin{proposition}
The equation $xy^2=y^2x$ does not transfer.
\begin{proof}
We construct a group in which the equation $xy^2=y^2x$ holds in a coset, but does not hold in any subgroup. Consider the finitely presented group $G$ on four generators $Z,Y,X,F$ given by the following equations.

~
\begin{enumerate}
    \item $F^2 = 1$,
    \item $XF=FX;\ XZ=ZX;\  XY=YX$,
    \item $ZF=FZ^{-1};\  YF=FY^{-1}$
    \item $YZ=ZYX$.
\end{enumerate}
The reader should check that the identities $Y^{-1}Z=ZY^{-1}X^{-1}$, $Y^{-1}Z^{-1}=Z^{-1}Y^{-1}x$ and $YZ^{-1}=Z^{-1}YX^{-1}$ also hold in the group~$G$. The equations above ensure that every word in $G$ can be written in a normal form $Z^{k}Y^{\ell}X^{m}F^{n}$ where $k,\ell,m \in \mathbb{Z}$ and $n \in \{0,1\}$.
For assume that $w \in G$ can be written in such a normal form $w=Z^{k}Y^{\ell}X^{m}F^{n}$. It suffices to show that the products $Fw,Xw,X^{-1}w,Yw,Y^{-1}w,Zw,Z^{-1}w$ can all be written in the same form.
\begin{itemize}
\item \textbf{Cases} $Fw, Xw, X^{-1}w, Zw, Z^{-1}w$\textbf{:} Trivial.
\item \textbf{Case} $Yw$\textbf{:} If $k \geq 0$, use the second and fourth cluster of identities to write $Yw=Z^{k}Y^{\ell+1}X^{m+k}F^{n}$. If $k$ is negative, use the identity~$YZ^{-1}=Z^{-1}YX^{-1}$ instead.
\item \textbf{Case} $Y^{-1}w$\textbf{:} If $k \geq 0$, use the second group and the identity~$Y^{-1}Z=ZY^{-1}X^{-1}$ to conclude that $Y^{-1}w=Z^{k}Y^{\ell-1}X^{m-k}F^{n}$. If $k$ is negative, use the identity~$Y^{-1}Z^{-1}=Z^{-1}Y^{-1}X$ to conclude that $Y^{-1}w=Z^{k}Y^{\ell-1}X^{m+k}F^{n}$.
\end{itemize}
We see that $G$ is a finitely generated infinite group of polynomial growth, so by Gromov's theorem it is virtually nilpotent. As in the proof of Proposition~\ref{prop:xyyy-equals-yyyx-tight}, we write every element of $G$ by giving the coefficients of its normal form, e.g. $(1,2,1,1)$ denotes the word $ZY^2XF$.

In what follows let $K,L,M \in \mathbb{Z}$. The five equations below define the group operation on the normal forms:
\begin{align}
&( 0, 0, 0, 0) \cdot (k,\ell,m,n) = (k,\ell,m,n) &\\
&(K, 0, 0, 0) \cdot (k,\ell,m,n) = ( 0, 0, 0, 0) \cdot (K+k, \ell,m,n) &\\
&(K,L, 0, 0) \cdot (k,\ell,m,n) = (K, 0, 0, 0) \cdot (k,L+\ell,m+Lk,n) &\\
&(K,L,M, 0) \cdot (k,\ell,m,n) = (K,L, 0, 0) \cdot (k,\ell,M+m,n) &\\
&(K,L,M, 1) \cdot (k,\ell,m,n) = (K,L,M, 0) \cdot (-k,-\ell,m,1+n) &
\end{align}
where we perform addition in the last component modulo 2. Now, any group element $w \in F\langle Z,Y,X \rangle$ can be written in the normal form $(k, \ell, m, 1)$ for some integers $k, \ell, m$. Using the equations above, we have:
\begin{align*}
w^2 &= (k, \ell, m, 1)\cdot (k, \ell, m, 1) &\text{ by (5)} \\
&= (k, \ell, m, 0)\cdot (-k, -\ell, m, 0) &\text{ by (4)} \\
&= (k, \ell, 0, 0)\cdot (-k, -\ell, 2m, 0) &\text{ by (3)} \\
&= (k, 0, 0, 0)\cdot (-k, 0, 2m, 0) &\text{ by (2)} \\
&= (0, 0, 0, 0)\cdot (0, 0, 2m, 0) &\text{ by (1)} \\
&= (0, 0, 2m, 0). &\text{ } \\
\end{align*}
We have $w^2 = (0, 0, 2m, 0) = X^{2m}$, and $X$ belongs to the center of $G$, so $w^2$ commutes with everything. Since $w \in F\langle Z,Y,X \rangle$ was an arbitrary element, we conclude that $xy^2=y^2x$ holds for all $x,y$ in the coset. Since the subgroup $\langle Z,Y,X \rangle$ has index 2 in $G$, the degree of satisfiability of the equation $xy^2=y^2x$ is at least $\frac{1}{2}$ with respect to any good density on the group $G$.
Finally, consider an arbitrary subgroup $H<G$. Consider elements of $H$ with normal form $(k,\ell,m,0)$. There are two possibilities:
\begin{enumerate}
\item Every such element has $\frac{l}{k}=c$ for some fixed ratio~$c$. Then we can use cosets of the form $Z^nH$ or $Y^nH$ to show that $H$ does not have finite index.
\item There are two elements $w_1 = (k_1,\ell_1,m_1,0) \in H$ and $w_2 = (k_2,\ell_2,m_2,0) \in H$ such that $\frac{l_1}{k_1} \neq \frac{l_2}{k_2}$. We calculate with the normal forms using equations (1-5), and quickly get that
\begin{align*}
w_1w_2w_2 = (k_1+2k_2,\ell_1+2\ell_2,m_1+2m_2+\ell_2k_2 + 2\ell_1k_2,0),\\
w_2w_2w_1 = (k_1+2k_2,\ell_1+2\ell_2,m_1+2m_2+\ell_2k_2 + 2\ell_2k_1,0).
\end{align*}
By uniqueness of normal forms, $w_1w_2w_2 = w_2w_2w_1$ precisely if $2\ell_1k_2 = 2\ell_2k_1$. But $\frac{\ell_2}{k_2} \neq \frac{\ell_1}{k_1}$ and therefore $xy^2=y^2x$ is not a law in $H$.
\end{enumerate}
We conclude that $xy^2=y^2x$ is not a law in any finite index subgroup $H<G$. Therefore the equation $xy^2=y^2x$ has positive degree of satisfiability in $G$ without being a virtual law.
\end{proof}
\end{proposition}

\section{Discussion}

\begin{body}
To what extent does Definition~\ref{def:transfer} characterize the equations which admit a direct analogue of Theorem~\ref{thm:antolin-commutativity}? In light of Theorem~\ref{thm:tointon-big}, it seems plausible that the correspondence could be exact (this would mean that naive generalizations of~\cite{antolin-infinite-dc} fail \textit{maximally}). Alas, we were not successful in proving more general transfer results (positive or negative) than the ones presented above. In particular, solutions to the following questions elude us:
\begin{enumerate}
    \item Does the equation $[x_1,x_2,\dots,x_n]=1$ transfer for $n>2$? Note that it is already known to satisfy the conclusion of Proposition~\ref{prop:transfer-big}.
    \item Given $p>2$, how can we construct a group $G$ where the equation $xy^p = y^px$ has positive degree of satisfiability but is not a virtual law?
\end{enumerate}
\end{body}

\begin{body}
The results of Section~\ref{sec:results} are strikingly elementary, and can be stated in familiar, group-theoretic language. We must note however that the language of logic allows for a far more unified treatment. We strongly hint at this Definition~\ref{def:degree-of-satisfiability}, which we formulate in terms of equations in the first-order language of group theory instead of the more conventional ``\textit{words on $n$ letters}''. As long as we work with equations, the difference between the algebraic (words) and logical (first-order formulas) seems insignificant. However, even if our interest lies in investigating the degree of satisfiability of group-theoretic \textit{equations}, some ways of obtaining results about these appear to lead us through satisfiability for general first-order formulae. 
\end{body}

\begin{body}
As a matter of fact, we already take such detours in the proofs of Propositions~\ref{prop:xyy-equals-yyx-gap}, \ref{prop:xyyy-equals-yyyx-gap} and \ref{prop:xy-equals-yX-gap}, where our roundabout ways are carefully obscured by the algebraic language. We can define the degree of satisfiability of arbitrary first-order formulae $\varphi$ in $k$ free variables (even with parameters) exactly as we did for equations in Definition~\ref{def:degree-of-satisfiability}. For example, the degree of satisfiability of the formula $\forall y. xy=yx$ in one free variable is just the following quantity:
$$ \frac{ \left|\SetComp{x \in G}{\forall y \in G. xy=yx}\right| }{ \left|G\right| }. $$
One quickly recognizes that the numerator denotes the size of $Z(G)$. From then on, it's obvious that $\forall y. xy=yx$ has finite satisfiability gap, with $\varepsilon=\frac{1}{2}$. In fact, we used this satisfiability gap repeatedly in the proofs of our equational gap results (whenever we mention Lagrange's theorem)! One could dismiss this as a coincidence, saying that the degree of satisfiability of $\forall y \in G. xy=yx$ helps with these proofs \textit{because} of the group-theoretic importance of centers and centralizers. However, the proof of Proposition~\ref{prop:xy-equals-yX-gap} further involves bounds on the sizes of sets $CL_G(a)$, first-order definable in parameters over $G$, that appear to have precious little group-theoretic significance. Given that non-trivial degree of satisfiability results \textit{can} be obtained in non-equational cases (with useful implications in the equational case), future investigation should not limit itself to the equational case, but focus on the full class of first-order formulae in the language of group theory.
\end{body}

{\small \paragraph{Acknowledgements.} The author thanks Benjamin M. Bumpus and Motiejus Valiunas for many helpful comments and suggestions.}
\pagebreak
\bibliography{thesis}
\bibliographystyle{ieeetr}


\end{document}